\documentclass[review]{elsarticle}
\usepackage{lineno,hyperref}
\usepackage[latin5]{inputenc}
\usepackage{amsmath}
\usepackage{graphicx}
\usepackage{subfigure}
\usepackage{caption}
\usepackage{float}
\modulolinenumbers[5]
\journal{.}
\begin{document}
\begin{frontmatter}
\title{Simulations of transport in one dimension}
\author{Alper Korkmaz}
\ead{akorkmaz@karatekin.edu.tr}
\address{Çankırı Karatekin University, Department of Mathematics, Çankırı, Turkey.}
\begin{abstract}
In this study, two initial boundary value problems for one dimensional advection-dispersion equation are solved by differential quadrature method based on sine cardinal functions. Pure advection problem modeling transport of conservative pollutants and fade out problem are simulated successfully by the proposed method. The time integration of the space discretized system is accomplished by using various single step and multi step methods covering forward, modified and improved Euler methods, Runge-Kutta, explicit Adams-Bashforth and implicit Adams-Moulton predictor-corrector methods of different orders. The errors between analytical and numerical solutions for both cases are measured by the use of discrete maximum norm. The numerical results are compared with some earlier results obtained by various methods.
\end{abstract}

\begin{keyword}
  Advection-Dispersion equation \sep transport \sep pollution \sep Sinc functions \sep  Differential quadrature method.
\end{keyword}

\end{frontmatter}

\linenumbers

\section{Introduction}

\noindent
Many physical phenomena in real world are modeled by various linear partial differential equations. Having both advection and dispersion (diffusion) terms in the Advection-Dispersion Equation (ADE) makes it a useful model for problems in various fields. Isenberg and Gutfinger examined a thin film of incompressible liquid draining down a vertical wall, motion for the film and unsteady heat transfer within the film \cite{isenberg}. Water transport in soils and dispersion in rivers and estuaries are also two well known studies modeled by the ADE\cite{parlange,chatwin}. Various problems including different types of the equation are used to model for the transient problems associated with flow through wellbores, geothermal production with reinjection, thermal energy storage in porous formations, thermal, hot fluid injection and energy extraction techniques for oil recovery, miscible flooding, oil recovery from hot dry rocks\cite{jochen}. A type of one dimensional form is used to describe uptake and desorption of solute diffusion into porous soil aggregates, lithofragments in sediments and aquifer materials in the sorptive\cite{grathwohl}. Solute transport problem by groundwater flow through isotropic and homogenous aquifer is also modeled by the ADE\cite{istok}. In transport phenomena in food processing, one dimensional unsteady diffusion in an isotropic medium, isothermal process, and the moisture content on a dry basis are studied with a different kind of the ADE\cite{jorge}. 

\noindent
Consider the initial boundary value problem for one dimensional form of the ADE
\begin{equation}
\frac{\partial u(x,t)}{\partial t}+\nu \frac{\partial u(x,t)}{\partial x}-\lambda \frac{\partial ^{2}u(x,t)}{\partial x^{2}}=0 \label{ad}
\end{equation}
with initial condition
\begin{equation*}
\begin{aligned}
u(x,0)&=f(x)
\end{aligned}
\end{equation*}
and boundary conditions 
\begin{equation*}
\begin{aligned}
u(a,t)&=b_1(t), \\
u(b,t)&=b_2(t)
\end{aligned}
\end{equation*}
over a finite interval $[a,b]$. This problem models transport of the quantity $u(x,t)$ of heat, fluid or related substances moving along $x-$axis with a constant flow velocity $\nu$ and the dispersion(diffusion) coefficient $\lambda$ \cite{bear,stocker}. 

\noindent
So far, various numerical methods have been applied to the ADE. Dağ et al. developed the least square finite element algorithm based on low degree B-spline shape functions (FEMLSF and FEMQSF) to solve transport problem modeled by the ADE \cite{dag1}. Szymkiewicz also solved a model problem described by the ADE via the combination of the spline functions and finite elements\cite{szym1}. Kadalbajaoo and Arora constructed a Taylor-Galerkin B-spline finite element algorithm to solve various initial boundary value problems for one dimensional advection-dispersion equation\cite{mohan1}. 

\noindent
Noye and Tan obtained the numerical solutions of the ADE by the third-order semi explicit finite difference method\cite{noye1}. Various two-level explicit and implicit finite difference methods covering the upwind explicit, the Lax-Wendroff, the modified Siemieniuch-Gladwell and the fourth-order method have been compared with each other on the numerical solutions of model problems including the ADE. Karahan solved various initial boundary value problems for the ADE by the use of implicit, third-order upwind and explicit finite difference methods\cite{karahan1,karahan2,karahan3}. Guraslan et al. developed a sixth-order compact finite difference method (CD6) combined with the fourth order Runge-Kutta method for numerical solution of three dynamic model problems\cite{guraslan1}. 

\noindent
Irk et al. set up a collocation method based on extended cubic B-spline functions (EXCBS)\cite{irk1}. In that study, pollutant transport through a channel problems modeled by the ADE with mixed boundary conditions were studied. Kaya developed a polynomial based differential quadrature algorithm to obtain numerical solutions of two initial boundary value problems including flood propagation in an open channel\cite{kaya1}. He also compared the obtained results with the explicit and implicit finite difference results. One more differential quadrature technique based on cubic B-spline functions (CSDQM) was developed for transport of conserved contaminant and fadeout problems in one dimension\cite{korkmaz1}. 

\noindent
Aim of this study is to obtain the numerical solutions of initial boundary value problems for the ADE in one dimension by differential quadrature method based on Sinc functions. The ordinary differential equation system obtained by the reduction of the ADE by differential quadrature method will be integrated for time variable by using various methods covering forward Euler(FORE) , improved polygon (modified Euler) method (IMPOLY), Heun (improved Euler) method(HEUN), classical Runge-Kutta methods of order two to four(RK2,RK3,RK4), implicit Rosenbrock method of third-fourth order(RB34), Gear single step method with Burlirsch-Stoer rational extrapolation(GB), FehlBerg Runge-Kutta method of order fifth order(RKF45), Runge-Kutta method with Cash-Karp coefficients of order four-five(RKCK45), Adams-Bashforth (AB4) and Adams-Moulton methods of order four(AM4). The first three initial steps of the iterations of AB4 and AM4 methods are calculated by RK4. In the predictor-corrector method AM4, the predictor method is chosen as AB4. 
\section{Numerical Method}
The Sinc functions
\begin{equation}
S_{m}(x)=\left\{ 
							\begin{array}{lcc}
										\dfrac{\sin{([\dfrac{x-m\Delta x}{\Delta x}]\pi)}}{[\dfrac{x-m\Delta x}{\Delta x}]\pi} & , & x \neq m\Delta x \\ 
										1 & , & x=m\Delta x \\ 
							\end{array}%
					\right.  \label{sinc}
\end{equation}
form a basis on the real line where $\Delta x$ is the equal node size, and $m$ is an integer \cite{stenger,carlson1,secer1,dehghan2}. The nodal values of sinc functions are described in \cite{dehghan2} as:
\begin{equation}
S_m(x_j)=\delta_{mj} %= \left\{ 							\begin{array}{lcc}										1& , & m=j\\ 										0 & , & m \neq j \\ 							\end{array}%					\right.  
\label{sincnodal}
\end{equation}

\noindent
Consider the series
\begin{equation}
C(u)(x)= \sum_{m=-\infty}^{\infty}u(m\Delta x)S_m(x)        
\label{cardinal}
\end{equation} 
for the function $u$ defined on ($-\infty ,\infty$). The function $C(u)(x)$ is named the cardinal of $u$ if it converges\cite{lund1}.
First two derivatives of Sinc function $S_m(x)$ are calculated as:
\begin{equation}
S_{m}^{\prime}(x)=\left\{ 
							\begin{array}{lcc}
										\dfrac{\dfrac{\pi}{\Delta x}(x-m\Delta x)\cos{\dfrac{x-m\Delta x}{\Delta x} \pi}-\sin{\dfrac{x-m\Delta x}{\Delta x} \pi} }{\dfrac{\pi}{\Delta x}(x-m\Delta x)^2}& , & x \neq m\Delta x \\ 
										0 & , & x=m\Delta x \\ 
							\end{array}%
					\right.  \label{sincd}
\end{equation}
\begin{equation}
S_{m}^{\prime \prime}(x)=\left\{ 
							\begin{array}{lcc}
							     \dfrac{-\dfrac{\pi}{\Delta x}\sin{\dfrac{x-m\Delta x}{\Delta x} \pi}}{x-m\Delta x}-\dfrac{2\cos{\dfrac{x-m\Delta x}{\Delta x} \pi}}{(x-m\Delta x)^2}+\dfrac{2\sin{\dfrac{x-m\Delta x}{\Delta x} \pi}}{\dfrac{\pi}{\Delta x}(x-m\Delta x)^3}& , & x \neq m\Delta x \\ 
									-\dfrac{{\pi}^2}{3\Delta x^2} & , & x=m\Delta x \\ 
							\end{array}%
					\right.  \label{sincdd}
\end{equation}

\noindent
Differential quadrature method (DQM) is a derivative approximation technique described as "\textit{the $p.$th order derivative of a function $u(x)$ at $x_m$ is approximated by finite weighted sum of nodal function values, i.e.,}
\begin{equation}
\left.  \frac{\partial u^{(p)}(x)}{\partial x^{(p)}}\right |_{x=x_{m}}=\sum \limits_{i=1}^{N}w_{mi}^{(p)}u(x_{i}),\quad m=1,2,\ldots ,N,  \label{Funda}
\end{equation}%
\textit{where the partion of the finite problem interval $[a,b]$ is $x_m=a+(m-1)\Delta x,m=1,2, \ldots , N$, $w_{mi}^{(p)}$ are the weights of nodal functional values for the $p.$ th order derivative approximation}\cite{bellman1}". The weights $w_{mi}^{(p)}$ are calculated using basis functions spanning the problem interval.
\subsection{Determination of the first order approximation weights}

\noindent
Letting $p=1$ in the fundamental differential quadrature derivative equation will lead to produce the weights of the first order derivative $w_{mi}^{(1)}$. The Sinc functions set $\{ S_m(x)\} _{m=1}^{m=N}$ forms a basis for the functions defined on $[x_1=a,b=x_N]$. In order to calculate the weights $w_{1i}^{(1)}$ of the node $x_1$, we substitute each Sinc basis functions into the fundamental differential quadrature equation \ref{Funda}. 
Substitution of $S_{1}(x)$ and using the functional and derivative values of it which can determined by using (\ref{sincd}) and (\ref{sincdd}) will lead the equation
\begin{equation}
\begin{aligned}
S_1^{\prime}(x_1)&=\sum \limits_{i=1}^{N}w_{1i}^{(1)}S_1(x_{i}) \\
     						 &=w_{11}^{(1)}{S_1(x_{1})}+w_{12}^{(1)}{S_1(x_{2})}+\ldots +w_{1N}^{(1)}{S_1(x_{N})} \\
								 &=w_{11}^{(1)}\delta_{11}+w_{12}^{(1)}\delta_{12}+\ldots +w_{1N}^{(1)}\delta_{1N} \\
							0	 &=w_{11}^{(1)} 		 
\end{aligned}
\label{w11}
\end{equation}
and will generate the weight $w_{11}^{(1)}$.
The weight $w_{12}^{(1)}$ can be calculated by substitution of $S_{2}(x)$ into Eq.(\ref{Funda}) as
\begin{equation}
\begin{aligned}
S_2^{\prime}(x_1)&=\sum \limits_{i=1}^{N}w_{2i}^{(1)}S_2(x_{i}) \\
     						 &=w_{11}^{(1)}{S_2(x_{1})}+w_{12}^{(1)}{S_2(x_{2})}+\ldots +w_{1N}^{(1)}{S_2(x_{N})} \\
								 &=w_{11}^{(1)}\delta_{21}+w_{12}^{(1)}\delta_{22}+\ldots +w_{1N}^{(1)}\delta_{2N} \\
							 \dfrac{(-1)^{2+1}}{\Delta x (1-2)} &=w_{12}^{(1)} 		 
\end{aligned}
\label{w12}
\end{equation}

\noindent
It can be concluded that the weights $w_{1i}^{(1)}$ focused on the first node $x_1$ can be determined by substitution of each Sinc functions $S_m(x),m=1,2, \ldots , N$ into the fundamental differential quadrature equation (\ref{Funda}) as
\begin{equation}
w_{1i}^{(1)}=\dfrac{(-1)^{1-i}}{\Delta x(1-i)},1\neq i  \label{a1i}
\end{equation}%
\begin{equation}
w_{11}^{(1)}=0  \label{a11}
\end{equation}%
When the weight $w_{mi}^{(1)}$ focussed on the node $x_m$ is wanted to be calculated, a general explicit formulation to determine it can be given as\cite{korkmaz2}:
\begin{equation}
w_{mi}^{(1)}=\dfrac{(-1)^{m-i}}{\Delta x(m-i)},m\neq i  \label{aij}
\end{equation}%
\begin{equation}
w_{mm}^{(1)}=0  \label{aii}
\end{equation}%
\subsection{Determination of the second order approximation weights}
Assuming $p=2$ and $m=1$ in the Eq.(\ref{Funda}) and using functional and derivative values of $S_1(x)$ will generate the equation
\begin{equation}
\begin{aligned}
S_1^{\prime \prime}(x_1)&=\sum \limits_{i=1}^{N}w_{1i}^{(2)}S_1(x_{i}) \\
     						 &=w_{11}^{(2)}{S_1(x_{1})}+w_{12}^{(2)}{S_1(x_{2})}+\ldots +w_{1N}^{(2)}{S_1(x_{N})} \\
								 &=w_{11}^{(2)}\delta_{11}+w_{12}^{(2)}\delta_{12}+\ldots +w_{1N}^{(2)}\delta_{1N} \\
							   \dfrac{-{\pi}^2}{3\Delta x^2}&=w_{11}^{(2)} 		 
\end{aligned}
\label{b12}
\end{equation}
Substitution of $S_2(x)$ into the fundamental differential quadrature equation (\ref{Funda}) will lead the equation
\begin{equation}
\begin{aligned}
S_2^{\prime \prime}(x_1)&=\sum \limits_{i=1}^{N}w_{1i}^{(2)}S_2(x_{i}) \\
     						 &=w_{11}^{(2)}{S_2(x_{1})}+w_{12}^{(2)}{S_2(x_{2})}+\ldots +w_{1N}^{(2)}{S_2(x_{N})} \\
								 &=w_{11}^{(2)}\delta_{21}+w_{12}^{(2)}\delta_{22}+\ldots +w_{1N}^{(2)}\delta_{2N} \\
							   2\dfrac{(-1)^{(2+1+1)}}{(\Delta x)^2(1-2)^2}&=w_{12}^{(2)} 		 
\end{aligned}
\label{b22}
\end{equation}
and will generate the weight $w_{12}^{(2)}$. In a general case the weights $w_{mi}^{(2)}$ focussed on the node $x_m$ of the second order derivative approximation can be written in an explicit form
\begin{equation}
w_{mi}^{(2)}=\dfrac{2(-1)^{m-i+1}}{\Delta x^{2}(m-i)^{2}},m\neq i  \label{bij}
\end{equation}%
\begin{equation}
w_{mm}^{(2)}=-\dfrac{\pi ^{2}}{3\Delta x^{2}}  \label{bii}
\end{equation} 
\noindent
\section{Discretization of the ADE}
\noindent
Replacing the space derivative terms by their DQM approximations in ADE (\ref{ad}) leads to an ordinary differential equation system of the form
\begin{equation}
\begin{aligned}
\left. \frac{\partial u(x,t)}{\partial t} \right |_{x=x_{m}}=-\nu \sum_{i=1}^{N}{w^{(1)}_{mi}u(x_{i},t)}+\lambda \sum_{i=1}^{N}{w^{(2)}_{mi}u(x_{i},t)}
\end{aligned} \label{disc1}
\end{equation}
where $w_{mi}^{(1)}$ and $w_{mi}^{(2)}$ are the weights of each $u(x_{i},t)$ for the first two derivative approximations at the node $x_m$. Since the nodal values of the function $u(x,t)$ at $x_1$ and $x_N$ are boundary conditions at both ends of the problem interval, then (\ref{disc1}) can be rewritten as
\begin{equation}
\begin{aligned}
\left. \frac{\partial u(x,t)}{\partial t} \right |_{x=x_{m}}&=\left [-\nu +\lambda\right ] w^{(1)}_{m1}b_{1}(t)+\left [-\nu+\lambda\right ] w^{(2)}_{mN}b_{2}(t) \\
&+\sum_{i=2}^{N-1} \left [ -\nu {w^{(1)}_{mi}}+\lambda w^{(2)}_{mi} \right ] u(x_{i},t)
\end{aligned} \label{disc2} 
\end{equation}
The fully space discretized system (\ref{disc2}) is integrated with respect to the time variable $t$ by using the time integration methods.
\section{Problems}
\noindent
In the process of application of numerical methods, the error between the numerical and the analytical results should be measured  to check the accuracy and the validity of the method. The measure of the error also provides a chance to compare the related method with the other ones. In this study, the discrete maximum norm $~_{(\Delta x, \Delta t)}L_{\infty}$ is used to measure the error between the numerical and the analytical solutions. This norm is defined as; 
$$
~_{(\Delta x, \Delta t)}L_{\infty}=\max_{2 \leq m \leq N-1}\left | {u^{\mathrm{a}}(x_{m},t)-u^{\mathrm{n}}(x_{m},t)}\right |
$$
where $u^{\mathrm{a}}(x_{m},t)$ and $u^{\mathrm{n}}(x_{m},t)$ are the analytical and the numerical solutions, respectively, at the node $x_m$ at a fixed time $t$ for the space and time step size $\Delta x$ and $\Delta t$.

\subsection{Transport with only Advection}
\begin{figure}[H]
    \subfigure[Simulation of the transport]{
   \includegraphics[scale =0.5] {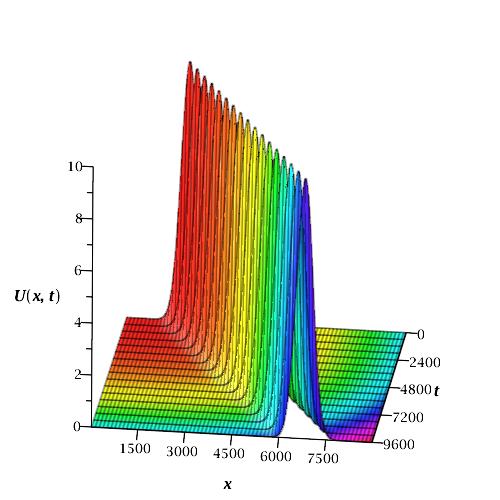}
   \label{fig:transportsim}
 }
 \subfigure[Maximum error as time goes with $\Delta x=25$ and $\Delta t=10$]{
   \includegraphics[scale =0.5] {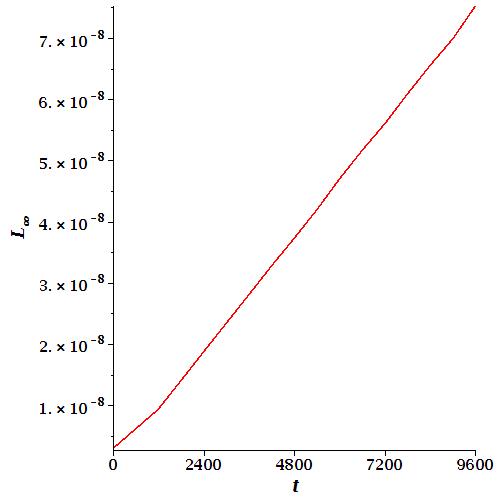}
   \label{fig:transporterror}
 }
 \caption{Transport of the initial quantity and the maximum errors during simulation}
\end{figure}

\noindent
The model problem for transport of a quantity of concentration along a channel is described as a pure advection initial boundary value problem for the ADE. The initial condition for the problem is derived by substituting $t=0$ into the analytical solution  
\begin{equation}
u(x,t)=10\exp (-\frac{1}{2\rho ^{2}}(x-\tilde{x}-\nu t)^{2})
\end{equation}
where $\rho$ and $\tilde{x}$ denote the standard deviation and the initial peak position of the bell-shaped quantity of 10 units height, respectively\cite{szym1,guraslan1,korkmaz1,dag1}. The solution represents motion of the initial quantity to the right along the channel of length 9 kilometers with a constant speed $\nu$. For the sake of comparison with the results stated in some earlier studies, the standard deviation $\rho=264$, the flow velocity $\nu=0.5\, m/s$ and the initial peak position $\tilde{x}=2$ referring the $2$ kilometers away from the left end of the channel are used as parameters to simulate the solutions. This choice of parameters moves the peak position of the initial quantity to $6.8$ kilometers far away from the left end of the channel at the simulation terminating time $9600$ seconds. The boundary conditions at both ends are selected as homogenous Dirichlet conditions over the problem interval $[0,9]$. The simulation of the transport obtained by SDQM-RKF45 with the parameters $\Delta x=25$ and $\Delta t=10$ is graphed in Fig  \ref{fig:transportsim}. The maximum error obtained by SDQM-RKF45 with the same parameters at some specific times are also recorded and depicted in Fig \ref{fig:transporterror}.
\noindent
Over the long simulation time, the solutions obtained by the SQDM seem stable and are in good agreement with the analytical ones. A comparison of the present results with the results in some earlier studies for various mesh sizes is tabulated in Table \ref{table:t1}. 

\noindent
When $\Delta x=200$ and $\Delta t=50$, the error obtained by  the SDQM-FORE is too high. The maximum norms are $1.15$ and $1.35$ for the CSDQM and the FEMLSF, respectively with the same parameters.  The results obtained by the methods SDQM-IMPOLY, SDQM-HEUN and SDQM-RK2 methods are accurate to one decimal digit like the results of the FEMQSF, the CD6 and the EXCBS. The RK3 and the AB4 methods have two decimal digits accuracy. The SDQM-RK4, the SDQM-GB, the SDQM-RKF45, the SDQM-RKCK45 and the SDQM-AB4 generate three decimal digit accurate results.

\noindent
The choice of $\Delta x$ and $\Delta t$ as $50$ causes to fail the low order the SDQM-FORE, the SDQM-IMPOLY, the SDQM-HEUN and the SDQM-RK2 and multi-step methods the AB4 and the AM4. The FEMLSF and the FEMQSF generates one decimal digit accurate results as the SDQM-RK3 has two decimal digits accuracy. The methods with three decimal digit accurate can be listed as the CSDQM and the EXCBS. The results obtained by the method CD6 are accurate to four decimal digits, the SDQM-RK4, the SDQM-RB34 and SDQM-RKF45 five decimal digits and the RKCK45 six decimal digits. The most accurate results obtained by the method SDQM-GB as eight decimal digits in this case. 

\noindent
Most of the methods applied for the time integration in this study, covering classical Runge-Kutta methods of order one to four, variations of Euler method and multi step methods, failed when $\Delta x$ is reduced to $25$ with fixed $Delta t=50$. The results obtained by the FEMLSF and the FEMQSF are accurate to one decimal digit, the CSDQM three decimal digits, and the CD6 four decimal digits. The accuracy of the results of SDQM-RB34 and the SDQM-RKF45 are five decimal digits as the best results again are obtained by the SDQM-GB as seven decimal digits accuracy.

\noindent
In the case reduction $\Delta t$ to $10$ with $\Delta x=25$, the SDQM-FORE method fails. The SDQM-IMPOLY, the SDQM-HEUN and the SDQM-RK2 generate two decimal digits accuracy as the results obtained by the SDQM-RK3 are accurate to four decimal digits, the SDQM-AB4 five decimal digits, the SDQM-RK4 and the SDQM-AM4 six decimal digits. The accuracy of the methods SDQM-RB34 and SDQM-GB are measured in seven decimal digits. The most accurate results for those parameters are produced by the methods SDQM-RK45 and SDQM-RKCK45 to eight decimal digits. Since the better results are obtained by the use of those parameters when compared with the results by EXCBS with $\Delta x= \Delta t=10$, we do not reduce the step sizes more.  

\begin{table}[H]
\scriptsize
\caption{Comparison of present results with some earlier ones for pure advection transport}
\begin{tabular}{lrrrrr}
\hline \hline 
   Method  &  $~_{(200,50)}L_{\infty}$ &  $~_{(50,50)}L_{\infty}$ &  $~_{(25,50)}L_{\infty}$ &  $~_{(25,10)}L_{\infty}$ &  $~_{(10,10)}L_{\infty}$\\
		\hline
  \hline
  SDQM-FORE   & 533.5714   &$\infty $&$\infty $ &$\infty $\\
  SDQM-IMPOLY  & 3.9486$\times 10^{-1}$  & $\infty $ &$\infty $ &1.7442$\times 10^{-2}$ \\
  SDQM-HEUN   & 3.9486$\times 10^{-1}$ & $\infty $&$\infty $ &1.5005$\times 10^{-2}$\\
  SDQM-RK2   & 3.9486$\times 10^{-1}$ & $\infty $ &$\infty$ &1.7442$\times 10^{-2}$ \\
  SDQM-RK3   & 1.9080$\times 10^{-2}$  & 1.8821$\times 10^{-2}$&$\infty $&1.5429$\times 10^{-4}$\\
  SDQM-RK4   & 1.9151$\times 10^{-3}$  & 7.0186$\times 10^{-5}$&$\infty $ &1.1436$\times 10^{-6}$\\
	SDQM-RB34   & 1.9182$\times 10^{-3}$  & 6.1214$\times 10^{-5}$ & 6.1275$\times 10^{-5}$&1.1967$\times 10^{-7}$ \\
	SDQM-GB   & 1.9183$\times 10^{-3}$  & 8.7642$\times 10^{-8}$ & 2.0875$\times 10^{-7}$& 1.1584$\times 10^{-7}$ \\
	SDQM-RKF45   & 1.9186$\times 10^{-3}$  & 1.8497$\times 10^{-5}$ & 1.8834$\times 10^{-5}$ & 7.5235$\times 10^{-8}$\\
	SDQM-RKCK45   & 1.9183$\times 10^{-3}$  & 3.0192$\times 10^{-6}$&23025.3677& 7.4091$\times 10^{-8}$\\
	SDQM-AB4  & 2.8709$\times 10^{-2}$  &  $\infty $ &  $\infty $& 4.6886$\times 10^{-5}$\\
	SDQM-AM4  & 2.5487$\times 10^{-3}$  &  $\infty $&  $\infty $& 3.5583$\times 10^{-6}$\\
	CSDQM \cite{korkmaz1}  & 1.15 &   8.00$\times 10^{-3}$& 1.00$\times 10^{-3}$\\
  FEMLSF \cite{dag1}  & 1.35 &   3.80 $\times 10^{-1}$& 3.77$\times 10^{-1}$\\
  FEMQSF \cite{dag1}  & 5.18 $\times 10^{-1}$&   3.73 $\times 10^{-1}$& 3.79$\times 10^{-1}$\\  
  CD6 \cite{guraslan1}  & 4.29$\times 10^{-1}$ & 8.00$\times 10^{-4}$ &7.00$\times 10^{-4}$ \\  
	EXCBS \cite{irk1}&6.07$\times 10^{-1}$&  2.20$\times 10^{-3}$&&&3.44$\times 10^{-6}$\\
\hline\hline
\end{tabular}
\label{table:t1}
\end{table}

\subsection{Transport with both Advection and Dispersion}

\begin{figure}[H]
 
   \subfigure[Fadeout of quantity as time goes]{
   \includegraphics[scale =0.5] {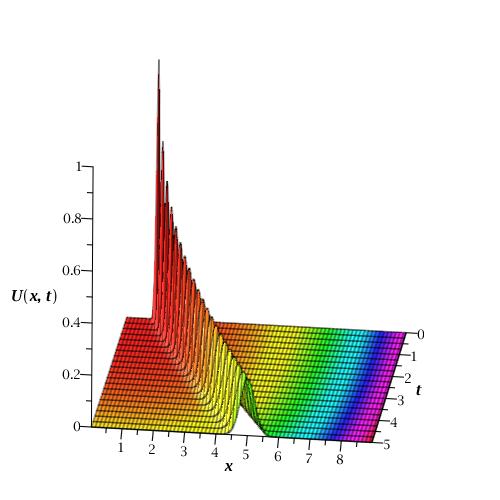}
   \label{fig:fadeout}
 }
 \subfigure[Maximum error as time goes]{
   \includegraphics[scale =0.5] {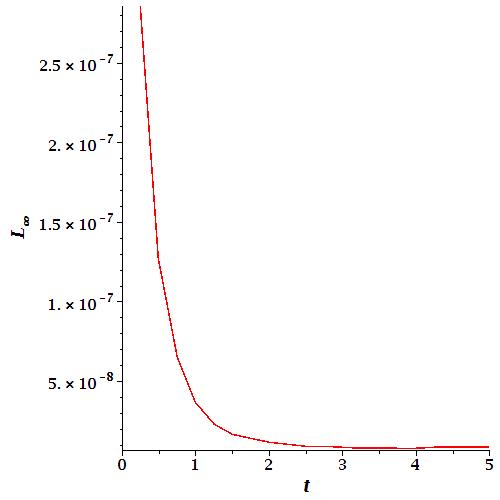}
   \label{fig:ferrors}
 }
 \caption{The fadeout of an initial quantity and the error at $t=5$}
\end{figure}

\noindent
The initial boundary value problem, constructed using both advection and dispersion terms together, models the fadeout of an initially solitary wave-shaped quantity while moving to the right along a channel as time goes. The analytical solution of this problem is given as
\begin{equation}
u(x,t)=\frac{1}{\sqrt{4t+1}}\exp \left( -\frac{(x-\tilde{x}-\nu t)^{2}}{%
\lambda (4t+1)}\right)
\end{equation}%
where $\tilde{x}$ is the initial peak position of the quantity of unit height moving with a constant velocity $\nu$\cite{noye1,mohan1}. The initial condition is chosen as 
\begin{equation}
u(x,0)=\exp \left( -\frac{(x-\tilde{x})^{2}}{\lambda}\right)
\end{equation}

\noindent
which can be obtained by substitution of $t=0$ into the analytical solution. The simulation is accomplished by assuming homogenous Dirichlet boundary conditions at both ends of the channel of length $9$ kilometers. The algorithm to simulate the solution of the problem is run up to the time $t=5$ seconds with the dispersion coefficient $\lambda=0.005$, the transport velocity $\nu=0.8 m/s$ and the initial peak position of the quantity $\tilde{x}=1$. The simulation of the motion and the maximum error-time graph are depicted in Fig \ref{fig:fadeout} and in Fig \ref{fig:ferrors}, respectively. The peak of the quantity reaches the fifth kilometers of the channel at the end of the simulation. This situation corresponds to the theoretical aspects of the solution owing to the value of $\nu$.

\noindent
A comparison of the results obtained by SQDM methods with the ones from the CSDQM method is also summarized for some various mesh sizes and fixed $\Delta t=0.0125$, Table \ref{table:t2}. When $\Delta x=0.2$, the results of all methods given in the table are as accurate as each other, namely to one decimal digit.

\noindent
When the mesh size is chosen as $0.1$, the results obtained from SDQM-FORE are one decimal digit accurate. This choice of $\Delta x$ causes two decimal digits accuracy for the method CSDQM(Method II). The results obtained by the CSDQM(Method I)has three decimal digits accuracy like all SDQM methods except SDQM-FORE. 

\noindent
In the case $\Delta x=0.05$, the results of SDQM-FORE has one decimal digit accuracy as the SDQM-AB4 fails. The accuracy of the results of the SDQM-IMPOLY, the SDQM-HEUN, the SDQM-RK2 and the CSDQM(Method I) are to three decimal digits. The methods SDQM-RK3, SDQM-RK4, SDQM-RB34, SDQM-GB, SDQM-RKF45, SDQM-RKCK45, SDQM-AM4, and CSDQM(Method II) have four decimal digits accurate.

\noindent
In the last case, we choose $\Delta x$ as $0.0025$. This choice of $\Delta x$ causes the methods SDQM-FORE, SDQM-RB34 and multi step methods to fail. The results obtained by the SDQM-IMPOLY, the SDQM-HEUN, and the SDQM-RK2 are accurate to three decimal digits accurate results, the CSDQM(Method I) four decimal digits, the SDQM-RK3 and the CSDQM(Method II) five decimal digits, the SDQM-RK4 seven decimal digits and the SDQM-GB and the SDQM-RKF45 eight decimal digits. The most accurate results are obtained by the SDQM-RKCK45 as nine decimal digits for this case.

\begin{table}[H]
\scriptsize
\caption{Comparison of the results with some earlier studies on the maximum error at $t=5$ for the fadeout problem}
\begin{tabular}{lrrrr}
\hline \hline 
   Method  &  $~_{(0.2,0.0125)}L_{\infty}$ &  $~_{(0.1,0.0125)}L_{\infty}$ &  $~_{(0.05,0.0125)}L_{\infty}$ &  $~_{(0.025,0.0125)}L_{\infty}$ \\
		\hline
  \hline
  SDQM-FORE   & 4.7876$\times 10^{-1}$ &2.2734$\times 10^{-1}$&2.2243$\times 10^{-1}$&$\infty $\\
  SDQM-IMPOLY  & 1.3818$\times 10^{-1}$ & 9.9836$\times 10^{-3}$ &1.6755$\times 10^{-3}$&1.6842$\times 10^{-3}$ \\
  SDQM-HEUN  & 1.3818$\times 10^{-1}$ & 9.9836$\times 10^{-3}$ &1.6755$\times 10^{-3}$&1.6842$\times 10^{-3}$ \\
  SDQM-RK2   & 1.3855$\times 10^{-1}$ & 9.9836$\times 10^{-3}$ &1.7655$\times 10^{-3}$&1.6842$\times 10^{-3}$ \\
  SDQM-RK3   & 1.3848$\times 10^{-1}$ & 9.9843$\times 10^{-3}$ &1.1087$\times 10^{-4}$&3.9909$\times 10^{-5}$ \\
  SDQM-RK4    & 1.3855$\times 10^{-1}$ & 9.9863$\times 10^{-3}$ &1.1070$\times 10^{-4}$&8.8121$\times 10^{-7}$ \\
  SDQM-RB34   & 1.3855$\times 10^{-1}$ & 9.9863$\times 10^{-3}$ &1.1071$\times 10^{-4}$&$\infty $ \\
  SDQM-GB   & 1.3855$\times 10^{-1}$ & 9.9863$\times 10^{-3}$ &1.1071$\times 10^{-4}$& 1.9130$\times 10^{-8}$ \\
  SDQM-RKF45   & 1.3855$\times 10^{-1}$ & 9.9863$\times 10^{-3}$ &1.1071$\times 10^{-4}$& 1.1869$\times 10^{-8}$ \\
	SDQM-RKCK45  & 1.3855$\times 10^{-1}$ & 9.9863$\times 10^{-3}$ &1.1071$\times 10^{-4}$& 8.6012$\times 10^{-9}$ \\
	SDQM-AB4  & 1.3856$\times 10^{-1}$ & 9.9860$\times 10^{-3}$ &$\infty $& $\infty $ \\
	SDQM-AM4  & 1.3855$\times 10^{-1}$ & 9.9864$\times 10^{-3}$ &1.1073$\times 10^{-4}$& $\infty $ \\
  CSDQM(Method I) \cite{korkmaz1}&  1.25$\times 10^{-1}$ & 6.95 $\times 10^{-3}$& 1.21$\times 10^{-3}$& 3.07$\times 10^{-4}$\\
  CSDQM(Method II) \cite{korkmaz1}& 1.36$\times 10^{-1}$& 1.45$\times 10^{-2}$ & 2.88$\times 10^{-4}$& 1.81$\times 10^{-5}$\\
 
\hline\hline
\end{tabular}
\label{table:t2}
\end{table}
\section{Conclusion}

\noindent
In the study, differential quadrature method based on sine cardinal functions is setup to solve the advection-dispersion equation numerically. The weight coefficients required for differential quadrature derivative approximations are computed in an explicit form. After discretization of the ADE in space by the DQM, and application of boundary conditions, the resultant ordinary differential equation system is integrated with respect to the time variable $t$ using various methods covering single step methods of different orders, and explicit Adams-Bahsforth and implicit Adams-Moulton multistep methods of order four. 
\noindent
In order to show the validity and accuracy of the numerical results, two initial boundary values problem are studied. The simulations and error distributions at the terminating times for both problems are depicted. The discrete maximum error norms measuring the error between the numerical and the analytical solutions are computed for various mesh and time step sizes. 
\noindent
A comparison of the results with each other and some results from different studies in literature is also accomplished by the comparison of norms. Comparisons also show that Sinc differential quadrature method generates acceptable, accurate and valid, better for some cases, solutions like the earlier solutions obtained by various methods in literature.


\section*{References}

\begin{thebibliography}{10}
\expandafter\ifx\csname url\endcsname\relax
  \def\url#1{\texttt{#1}}\fi
\expandafter\ifx\csname urlprefix\endcsname\relax\def\urlprefix{URL }\fi
\expandafter\ifx\csname href\endcsname\relax
  \def\href#1#2{#2} \def\path#1{#1}\fi

\bibitem{isenberg}
J.~Isenberg, C.~Gutfinger, Heat transfer to a draining film, International
  Journal of Heat and Mass Transfer 16 (1972) 505--512.

\bibitem{parlange}
J.~Y. Parlange, Water transport in soils, Annual Review of Fluid Mechanics 19
  (1980) 77--102.

\bibitem{chatwin}
P.~Chatwin, C.~Allen, Mathematical models of dispersion in rivers and
  estuaries, Annual Review of Fluid Mechanics 17 (1985) 119--149.

\bibitem{jochen}
J.~Bundschuh, M.~C.~S. Arriaga, Introduction to the Numerical Modeling of
  Groundwater and Geothermal Systems: Fundamentals of Mass, Energy and Solute
  Transport in Poroelastic Rocks, CRC Press, Boca Raton, 2010.

\bibitem{grathwohl}
P.~Grathwohl, DIFFUSION IN NATURAL POROUS MEDIA:Contaminant Transport,
  Sorption/Desorption and Dissolution Kinetics, Kluwer Academic Publishers,
  Newyork, 1998.

\bibitem{istok}
J.~Istok, Groundwater Modelling by the Finite Element Method, American
  Geophysical Union, NW Washington DC, 1989.

\bibitem{jorge}
S.~M. Alzamora, A.~Nieto, M.~A. Castro, Structural effects of blanching and
  osmotic dehydration pretreatments on air drying kinetics of fruit tissues,
  in: J.~Welti-Chanes, J.~F. V{\'e}lez-Ruiz, G.~V. Barbosa-C{\'a}novas (Eds.),
  Transport Phenomena in Food Processing, CRC Press LLC, Boca Raton, 2003.

\bibitem{bear}
J.~Bear, A.~Verjuit, Modeling Groundwater Flow and Pollution, D. Reidel
  Publishing Company, Dordrecht, 1987.

\bibitem{stocker}
T.~Stocker, Introduction to Climate Modelling, Springer, Heidelberg, 2011.

\bibitem{dag1}
I.~Da\u{g}, D.~Irk, M.~Tombul, Least-squares finite element method for the
  advection-diffusion equation, Applied Mathematics and Computation 173 (2006)
  554--565.

\bibitem{szym1}
R.~Szymkiewicz, Solution of the advection-diffusion equation using the spline
  function and finite elements, Communications in Numerical Methods in
  Engineering 9 (1993) 197--206.

\bibitem{mohan1}
K.~K. Mohan, P.~Arora, Taylor-{G}alerkin {B}-spline finite element method for
  the one-dimensional advection-diffusion equation, Numerical Methods for
  Partial Differential Equations 26 (2010) 1206--1223.

\bibitem{noye1}
B.~J.~Noye and H. H.~Tan, A third-order semi-implicit finite difference method for
  solving the one-dimensional convection-diffusion equation, Communications in
  Numerical Methods in Engineering 26 (1988) 1615--1629.

\bibitem{karahan1}
H.~Karahan, Implicit finite difference techniques for the advection-diffusion
  equation using spreadsheets, Advances in Engineering Software 37 (2006)
  601--608.

\bibitem{karahan2}
H.~Karahan, A third-order upwind scheme for the advection-diffusion equation
  using spreadsheets, Advances in Engineering Software 38 (2006) 688--697.

\bibitem{karahan3}
H.~Karahan, Unconditionally stable explicit finite difference technique for the
  advection-diffusion equation using spreadsheets, Advances in Engineering
  Software 38 (2007) 80--86.

\bibitem{guraslan1}
G.~Gurarslan, H.~Karahan, D.~Alkaya, M.~Sari, M.~Yasar, Numerical solution of
  advection-diffusion equation using a sixth-order compact finite difference
  method, Mathematical Problems in Engineering Article ID 672936 (2013) 1--7.

\bibitem{irk1}
D.~Irk, I.~Da\u{g}, M.~Tombul, Extended cubic b-spline solution of the
  advection-diffusion equation, KSCE Journal of Civil Engineering 19 (2014)
  929--934.

\bibitem{kaya1}
B.~Kaya, Solution of the advection-diffusion equation using the differential
  quadrature method, KSCE Journal of Civil Engineering 14 (2009) 69--75.

\bibitem{korkmaz1}
A.~Korkmaz, I.~Da\u{g}, Cubic B-spline differential quadrature methods for the advection-diffusion equation, International Journal of Numerical Methods
  for Heat {\&} Fluid Flow 22 (2012) 1021--1036.

\bibitem{stenger}
F.~Stenger, Numerical Methods Based on Sinc and Analytic Functions, Springer,
  New York, 1993.

\bibitem{carlson1}
T.~S. Carlson, J.~Dockery, J.~Lund, A sinc-collocation method for initial
  boundary value problems, Mathematics of Computation 66 (1997) 215--235.

\bibitem{secer1}
A.~Secer, Numerical solution and simulation of second-order parabolic pdes with
  sinc-galerkin method using maple, Abstract and Applied Mathematics Article ID
  686483 (2013) 1--10.

\bibitem{dehghan2}
M.~Dehghan, A.~Saadatmandi, The numerical solution of a nonlinear system of
  second-order boundary value problems using the sinc-collocation method,
  Mathematical and Computer Modelling 46 (2007) 1434--1441.

\bibitem{lund1}
J.~Lund, K.~L. Bowers, Sinc Methods for Quadrature and Differential Equations,
  SIAM, Philadelphia, 1992.

\bibitem{bellman1}
R.~Bellman, B.~G. Kashef, J.~Casti, Differential quadrature: A tecnique for the
  rapid solution of nonlinear differential equations, Journal of Computational
  Physics 10 (1972) 40--52.

\bibitem{korkmaz2}
A.~Korkmaz, I.~Da\u{g}, Shock wave simulations using sinc differential
  quadrature method, Engineering Computations 28 (2011) 654 -- 674.

\end{thebibliography}
\end{document}